\documentclass[12pt]{article}   
\usepackage{amsfonts}   

\begin{document}

\newtheorem{lemma}{Lemma}[section]
\newtheorem{theo}[lemma]{Theorem}
\newtheorem{coro}[lemma]{Corollary}
\newtheorem{rema}[lemma]{Remark}
\newtheorem{propos}[lemma]{Proposition}

\newcommand{\chu}{{\cal M}}
\newcommand{\shu}{\mbox{$ S( H^{\infty} ) $}}  
\newcommand{\rr}{\mbox{$    \rightarrow   $}}
\newcommand{\papa}{ H^{\infty} }
\newcommand{\disc}{ {\Bbb D} }
\newcommand{\ov}{ \overline }
\newcommand{\zinf}{ Z_{\infty} }
\newcommand{\inter}{ \mbox{int\/} }
\newcommand{\diaend}{$\square$}
\newcommand{\om}{\omega}
\newcommand{\vsn}{\vspace{1mm}}
\newcommand{\ma}{K}
\newcommand{\hull}{\mbox{hull}\,}

\newcommand{\alphak}{\alpha_{k}}
\newcommand{\rhok}{\rho_{k}}
\newcommand{\deltak}{\delta_{k}}
\newcommand{\gammak}{\gamma_{k}}
\newcommand{\xik}{\xi_{k}}
\newcommand{\betak}{\beta_{k}}
\newcommand{\varj}{\varepsilon_{j}}
\newcommand{\vark}{\varepsilon_{k}}
\newcommand{\sigmak}{\sigma_{k}}
\newcommand{\vaci}{\mbox{ }  }
\newcommand{\lambdak}{\lambda_{k}}
\newcommand{\diamh}{\mbox{diam}_{h}} 
\newcommand{\sspace}{\,}
\newcommand{\dist}{\mbox{dist}} 
\newcommand{\vsa}{\vspace{1.4mm}} 

\newcommand{\labeth}[1]{\label{{#1}.th}}
\newcommand{\refth}[1]{\ref{{#1}.th}}
\newcommand{\laba}{\label}     

\hyphenation{cha-rac-te-ri-za-tion}  \hyphenation{res-pec-ti-ve-ly}
\hyphenation{en-cou-ra-ge-ment}    \hyphenation{re-gu-la-ri-zed}
\hyphenation{subs-tan-tia-lly}          \hyphenation{pa-ra-me-ters}  
\hyphenation{na-tu-ral}    \hyphenation{co-ro-lla-ry}      \hyphenation{o-ther-wise} 
\hyphenation{sa-tis-fies}  \hyphenation{po-si-ti-ve}
\hyphenation{pro-ducts}  \hyphenation{pro-blem}         \hyphenation{ve-ri-fy}
\hyphenation{ge-ne-ral}   \hyphenation{pa-ra-me-ter}   \hyphenation{boun-ded}
\hyphenation{pseudo-hyper-bo-llic}   \hyphenation{e-ve-ry} 
\hyphenation{iso-me-try}   \hyphenation{tri-vial}    \hyphenation{pro-duct}
\hyphenation{pro-ducts}

\renewcommand{\thesection}{\arabic{section}}
\renewcommand{\theequation}{\thesection.\arabic{equation}}
\newcommand{\equnew}{\setcounter{equation}{0}}

\newcommand{\calo}{  {\cal O} }
\newcommand{\im}{\mbox{Im} \, }
\newcommand{\papac}[1]{   H^{\infty}|_{ \mbox{{\scriptsize supp}}\,  {#1} }   }
\newcommand{\chuc}[1]{  M( H^{\infty}|_{  \mbox{{\scriptsize supp}}\,  {#1}   })   }
\newcommand{\supp}{ \mbox{supp}\,  } 
\newcommand{\Om}{\Omega} 
\newcommand{\zd}{Z_{D}} 
\newcommand{\calcu}{\Theta} 
\newcommand{\calg}{\cal G} 
\newcommand{\eiti}{e^{i\theta}} 
\newcommand{\ceiti}{e^{-i\theta}} 
\newcommand{\hypdi}{\Delta }   
\newcommand{\hh}{h}    

\title{On the analytic structure of  the $\papa$ maximal ideal space} 
\author{Daniel Su\'{a}rez} 

\date{\mbox{}}
\maketitle 

\begin{quotation}
\noindent
\mbox{ } \hfill      {\sc Abstract}  \hfill \mbox{ } \\  
\footnotetext{ 2000 Mathematics
Subject Classification: primary 46J15, secondary 30D50. Key words:
Gleason parts, analytic extension. }       \hfill \mbox{ }\\
{\small We characterize the algebra $\papa \circ L_{m}$, where 
$m$ is a point of the maximal ideal space of $\papa$  with nontrivial 
Gleason part $P(m)$ and $L_{m} : \disc \rr P(m)$ is the coordinate 
Hoffman map. 
In particular, it is shown that for any continuous function 
$f: P(m) \rr {\Bbb C}$ with $f\circ L_{m} \in \papa$ there 
exists $F\in \papa$ such that $F|_{P(m)} = f$.} \\ 
\end{quotation}
%
\section*{Introduction - Preliminaries}

\noindent        
The Gelfand transform represents $\papa$ as an algebra of continuous functions 
on its maximal ideal space $\chu$ 
(provided with the weak star topology) 
via the formula  
$\hat{f} (\varphi ) \stackrel{\mbox{\scriptsize{def}}}{=}  \varphi (f)$, 
where $f\in \papa$ and $\varphi \in \chu$. 
We will not write the hat of $f$ unless the contrary is stated. 
Beginning with a seminal paper of Hoffman \cite{hof}, 
many papers have studied the analytic behavior of 
$\papa$ on  parts of  $\chu$ others than the disk  
(see \cite{bud},$\,$\cite{gor},$\,$\cite{g-l-m},$\,$\cite{iz} and \cite{su4}).    
While the whole picture seems to be  unreachable, 
the present paper intends to throw some light into this 
never-ending program. A more precise statement of our result 
(Thm.$\,$\refth{big} and Coro.$\,$\refth{thm?}) 
will require to develop some notation and machinery.  

The pseudohyperbolic metric for $x, y \in \chu$ is defined by 
\[ 
\rho (x, y)  =   \sup \{  | f(y) |  :   f\in  \papa ,   \     \| f \| =1   \mbox{ and }    f(x) = 0  \} ,   
\] 
which for $z, \om \in \disc$ reduces to $\rho (z, \om ) = | z-\om | / |1-\ov{\om}z|$.    
The Gleason part of $m\in \chu$ is 
$P(m)  \stackrel{\mbox{\scriptsize{def}}}{=}  \{  x\in \chu :  \rho (m,x ) < 1 \}$. 
Clearly $\disc$ is a Gleason part. 
If $z_{0}\in \disc$, we can think of
the analytic function 
\[
L_{z_{0}}(z) =
\frac{z+z_{0}}{1+ \ov{z}_{0} z } ,  \ z\in \disc
\]  
%
as mapping $\disc$ into $\chu$. 
In \cite{hof} Hoffman proved that if $m\in \chu$ and 
$(z_{\alpha})$ is a net in $\disc$ converging to $m$,  
then the net $L_{z_{\alpha}}$ tends in the space $\chu^{\disc}$ 
(i.e., pointwise) 
to some analytic map $L_{m}$ from $\disc$  onto $P(m)$  
such that $L_{m}(0) = m$. Here `analytic' means 
that $f\circ L_{m} \in \papa$ for every $f\in \papa$. 
The map $L_{m}$ does not depend on the particular choice of  the 
net $(z_{\alpha})$ that tends to $m$. 
A Blaschke product $b$ with zero sequence $\{ z_{n} \}$ satisfying 
\[              \delta (b) =    
\delta ( \{ z_{n} \} )  \stackrel{\mbox{\scriptsize{def}}}{=}  
\inf_{n} \prod_{j \, : \, z_{j} \not = z_{n}}  \rho (z_{n} , z_{j}) > 0
\]  
is called an interpolating Blaschke product and $\{ z_{n} \}$ is called an 
interpolating sequence.  
Let $\calg$ denote the set of 
points in $\chu$ that lie in the closure of some interpolating sequence. 
If $m\in \chu \setminus \calg$ then $P(m) = \{ m \}$ and hence $L_{m}$ 
is a constant map. If $m\in \calg$ then $L_{m}$ is one-to-one, 
meaning that $P(m)$ is an analytic disk in $\chu$.  
Hoffman also realized 
that even when $P(m)$ is a disk, there are cases in which 
$L_{m}$ is a homeomorphism  and cases in which it is not. 

By an abstract version of  Schwarz's lemma \cite[p.$\,$162]{sto},  
any connected portion of  $\chu$ provided with a nontrivial  
analytic structure  must be contained in some $P(m)$ with $m\in \calg$. 
In order to understand the analytic structure of $\chu$ it is then fundamental 
to study the Hoffman algebras $\papa \circ L_{m}$, where $m\in \calg$. 

In \cite{su4} it is proved that  $\papa \circ L_{m}$ is a closed subalgebra 
of $\papa$ and that they coincide when $P(m)$ is a homeomorphic disk 
(i.e., $L_{m}$ is a homeomorphism). 
Particular versions of the last result were 
obtained in \cite[pp.$\,$106-107]{hof} and \cite[Coro.$\,$3.3]{g-l-m}.    
On the other hand, when $P(m)$ is not a homeomorphic disk it is well known 
that the identity function is not in $\papa \circ L_{m}$, meaning that this 
algebra is properly contained in $\papa$. But, what is it? 
We provide an answer to this question by giving several 
characterizations of  $\papa \circ L_{m}$, 
the most natural being 
\[  \papa \circ L_{m}  = 
\papa \cap [ C(P(m), {\Bbb C}) \circ L_m ] ,  
\]   
where $C(P(m), {\Bbb C})$ is the algebra of continuous maps from $P(m)$ 
into ${\Bbb C}$. 
The inclusion $\subseteq$ is trivial, but the proof of the other 
inclusion turned out to be very difficult.  
We can look at the equality as an extension result; it says that 
for every continuous function $f$ on $P(m)$ such that $f\circ L_{m} \in \papa$ 
there is an extension $F\in \papa$ of $f$ (i.e.,  $F|_{P(m)}  = f$). 
By the above comments, the result is new only for non-homeomorphic disks,  
but the argument here works in general. 
However, 
the technical complications 
introduced by considering non-homeomorphic disks make the proof 
much more difficult and longer  than in \cite{su4}. 


\section{Algebraic properties of Hoffman maps} 

\noindent 
Let $\tau : \disc \rr \chu$ be an analytic function. As before, this means that 
$f\circ \tau \in \papa$ for all $f\in \papa$.  
We can extend $\tau$ to a continuous map $\tau^{\ast}: \chu \rr \chu$ by the formula 
$\tau^{\ast} (\varphi )  (f)        \stackrel{\mbox{\scriptsize{def}}}{=} 
\varphi ( f\circ \tau )$, where $\varphi \in \chu$.   
Two particular cases will be of interest here. 
If  $\tau$ is an analytic self-map of $\disc$  
then we can think of $\tau$ as mapping $\disc$ into 
$\chu$ and consider its extension $\tau^{\ast}$. 
We can do this with any automorphism of $\disc$, which 
therefore induces a homeomorphism 
from $\chu$ onto $\chu$. In particular, if $\lambda \in \partial \disc$,  
the rotation $z \mapsto \lambda z$ ($z\in \disc$) extends to $\chu$ in this way. 
From now on, for $\varphi \in \chu$ we simply write $\lambda \varphi$ 
for this `$\lambda$-rotation' in $\chu$.  We  point out that even when 
each such rotation is a homeomorphism, the action of the group 
$\partial \disc$  into $\chu$ is not continuous \cite[pp.$\,$164-165]{hof.b}. 
The other relevant case for the paper is $L_{m}$ (for $m\in \calg$). 
The extension $L^{\ast}_{m}$ maps $\chu$ onto $\ov{P(m)}$, and 
$P(m)$ is a homeomorphic disk if and only if  $L^{\ast}_{m}$ 
is one-to-one \cite[Sect.$\,$3]{su4}. 
We will also denote this extension  by $L_{m}$, 
where the meaning will be clear from the context. 

The inclusion of the disk algebra in $\papa$ induces 
a natural projection $\pi : \chu \rr \ov{\disc}$. The fiber of a point 
$\om \in \partial \disc$ is $\pi^{-1} (\om ) \subset \chu$.  
Let $x, y \in \chu$ and let $(z_{\alpha})$ be a net in $\disc$ so that 
$y = \lim z_{\alpha}$. We claim that the limit of 
$(1+ \ov{\pi (x)} z_{\alpha})  /    (1+ \pi (x) \ov{z_{\alpha}})$ 
always exists (in $\ov{\disc}$) and it is independent of the net $(z_{\alpha})$. 
A rigorous statement would say that the above limit exists when $\pi (z_{\alpha})$ 
is in place of  $z_{\alpha}$; but since $\pi$ identifies $\disc$ with $\pi (\disc )$, 
no harm is done with the appropriate mind adjustment.  

It is clear that the limit exists whenever the denominator does not tend to zero. So, 
suppose that $|\pi (x)| =1$ and $\pi (y) = - \pi (x)$. 
Write $z_{\alpha} = - \pi (x) (1- r_{\alpha} e^{i\theta_{\alpha}})$, with 
$-\pi /2 < \theta_{\alpha} < \pi /2$  and  $0<r_{\alpha} \rr 0$.     
The point $y\in \chu$ is either a nontangential point, 
in which case $\theta_{\alpha} \rr \theta \in (-\pi /2 , \pi /2)$ 
(see \cite[pp.$\,$107-108]{hof}), 
or it is a tangential point and 
$\theta_{\alpha}$ accumulates in $\{  -\pi /2 , \      \pi /2  \}$.  
In both cases  
\[ 
\lim \frac{  1+ \ov{\pi (x)} z_{\alpha}  }{  1+ \pi (x) \ov{z_{\alpha}}  }   = 
\lim \frac{  1-(1- r_{\alpha} e^{i\theta_{\alpha}})    }{   
		 1-(1- r_{\alpha} e^{-i\theta_{\alpha}})  }    = 
\lim  e^{i2\theta_{\alpha}} ,
\]
proving our claim. \\

\noindent 
{\sc Definition.} Let $\lambda : \chu \times \chu  \,    \rr  \,   \partial \disc \,$ 
be the function 
\[   
\lambda (x,y) = 
\lim_{\alpha} \frac{  1+ \ov{\pi (x)} z_{\alpha}  }{  1+ \pi (x) \ov{z_{\alpha}}  } 
\ \ \    (x , y \in \chu ), 
\]    
where $(z_{\alpha})$ is any net in $\disc$ that tends to $y$. \\

\noindent 
Observe that if $x, y \in \chu$  do not satisfy the extreme conditions 
$|\pi (x)| =1$ and $\pi (x) = -\pi (y)$ then 
\[   \lambda (x,y) = 
\frac{  1+ \ov{\pi (x)} \pi (y)  }{  1+ \pi (x) \ov{\pi (y)}  } , 
\]
and this expression reduces to $\ov{\pi (x)} \pi (y)$ when $|\pi (x)| =1 = |\pi (y)|$. 
We will use indistinctly the notations $\lambda (x,y)$  or  $\lambda_{x,y}$ 
to denote this function. 
A word of warning: the function $\lambda$ was introduced 
by Budde in \cite{bud} for the purpose of proving the same 
result given in Proposition \refth{compo} below. 
However, the value of  $\lambda (x,y)$ stated in \cite{bud} 
when $|\pi (x)| =|\pi (y)| =1$ is $\ov{\pi (x)}\pi (y)$, 
therefore overlooking the pathological behavior of 
$\lambda$ when $\pi (x) = -\pi (y)$. Fortunately, all the proofs and results   
in \cite{bud} remain valid by only adjusting $\lambda$ to its right value. 

In \cite[Lemma 1.8]{g-l-m} Gorkin, Lingerberg and Mortini 
proved that if $m\in \calg$ and $b$ is an interpolating Blaschke product then 
$b\circ L_{m} = B f$, where $B$ is an interpolating Blaschke product 
and $f\in (\papa )^{-1}$. 
This fact will be used frequently along the paper.  
We will also need a result of Budde 
\cite{bud} stating that if $\varphi \in \chu$ has trivial Gleason part then 
$L_{m} (\varphi )$  also has trivial Gleason part. 
In symbols, $L_{m}^{-1} (\calg ) \subset \calg$.

\begin{propos}                          \labeth{compo}
Let $m, y \in \chu$.
Then  
\begin{equation}          \laba{pong} 
L_{m} \circ L_{y} (\lambda (m,y) z) =  L_{ L_{m}(y) }(z)  
 \ \ \mbox{ for all }\ \  z\in \disc . 
\end{equation} 
\end{propos} 
{\em Proof.} 
Suppose first that $m$ or $y$ (or both) is not in $\calg$. 
Then one of the maps $L_{m}$ or $L_{y}$ is constant and the left member of 
(\ref{pong}) is the constant map $L_{m}(y)$. If $y\not \in \calg$ then 
Budde's result asserts that $L_{m}(y)\not \in \calg$, which is trivially 
the case  if $m\not \in \calg$, too. 
Therefore also the right member of (\ref{pong}) 
is the constant map $L_{m}(y)$. 

Suppose now that $m, y\in \calg$ and let $\om  , \xi$ and $z$ in $\disc$. 
An elementary calculation shows that 
\begin{equation}    \laba{lom} 
L_{\om} ( L_{\xi} (z) ) = L_{ L_{\om}(\xi )} (\ov{\lambda (\om , \xi )} z) .  
\end{equation} 
Replace $\om$ in 
(\ref{lom}) by a net $(\om_{\alpha})$ in $\disc$ tending to $m$. 
Then the first member of (\ref{lom}) tends to $L_{m} (L_{\xi} (z))$ and  
$L_{\om_{\alpha}}(\xi ) \rr L_{m}(\xi )$. Consequently  
\begin{equation}             \laba{sub} 
 L_{  L_{\om_{\alpha}}(\xi )  }  \rr L_{  L_{m}(\xi )  } 
\ \ \mbox{ pointwise on }\ \ \disc .
\end{equation} 
It is clear that the constants 
$\lambda_{\alpha} =\lambda (\om_{\alpha} ,\xi )$ 
tend to $\lambda (m, \xi )$.   
Using that for $x\in \calg$ the map $L_{x}$ is an isometry on $\disc$ 
with respect to $\rho$ \cite[p.$\,$105]{hof} we get  
%
\begin{eqnarray}           \label{rhomu} 
\rho ( L_{  L_{\om_{\alpha}}(\xi )  }  (\ov{\lambda}_{\alpha} z)  \,  ,   \, 
	 L_{  L_{\om_{\alpha}}(\xi )  }  (\ov{\lambda (m, \xi )}   z)  )   
	 &  = & 
\rho (  \ov{\lambda}_{\alpha}z,  \ov{\lambda (m, \xi )} z  )  \nonumber \\
       &   \leq   & 
\frac{   |\ov{\lambda}_{\alpha}- \ov{\lambda (m, \xi )}|  }{ 1-|z|^{2} }  
\rr  0 .
\end{eqnarray}
Hence, the lower semicontinuity of $\rho$ 
(see \cite[Thm.$\,$6.2]{hof}) together with (\ref{lom}), 
(\ref{sub}) and (\ref{rhomu}) yields 
 \begin{equation}     \label{puchch} 
\rho (  L_{m } (L_{\xi} (z))  , L_{ L_{m}(\xi )} (\ov{\lambda (m, \xi )}   z)  ) 
\leq  
\lim  \hspace{0.15mm}    
\rho (  L_{  L_{\om_{\alpha}}(\xi )  }  (\ov{\lambda}_{\alpha} z)  , 
           L_{  L_{\om_{\alpha}}(\xi )  }  (\ov{\lambda (m, \xi )}  z)  )     
=0 .   
\end{equation}    
That is,  $  L_{m } (L_{\xi} (z))  = L_{ L_{m}(\xi )} ( \ov{\lambda (m, \xi )}  z)$ 
for every $\xi , z\in \disc$.   
Now replace $\xi$ by a net $(\xi_{\alpha}  )$ in $\disc$ tending to $y$. 
By the continuity of $L_{m}$ on $\chu$ then 
$L_{m} ( L_{ \xi_{\alpha} } (z)  ) \rr  L_{m}( L_{y}(z) )$ and 
$L_{m}(\xi_{\alpha})  \rr  L_{m} (y)$. 
Since the map $x\mapsto L_{x}$ is continuous from $\chu$ into $\chu^\disc$, then 
we also have 
\[         L_{ L_{m}(\xi_{\alpha}) } \rr  L_{ L_{m} (y)} 
\ \ \mbox{ pointwise on }\ \ \disc .
\] 
Since $L_{m} (\xi_{\alpha}) \in P(m) \subset \calg$ for every $\xi_{\alpha}$,  
then $L_{     L_{m}(\xi_{\alpha})    }$ are isometries with respect to $\rho$.  
In addition, $\lambda (m, \xi_{\alpha})  \rr  \lambda (m,y)$ by definition. 
Therefore the same argument as in  (\ref{rhomu}) and (\ref{puchch}) yields  
\[  L_{m} (L_{y} (z)) =  L_{ L_{m}(y) } (\ov{ \lambda (m,y) } z) .
\]
The proposition follows replacing $z$ by $\lambda (m,y)  z$.    \diaend  


\begin{coro}[Budde]      \labeth{o-o} 
Let $m\in {\cal G}$ and $\xi \in \chu$ such that $L_{m}(\xi ) \in P(m)$. 
Then $L_{m}$ maps $P(\xi )$  onto $P(m)$ in a one-to-one 
fashion. 
\end{coro}
{\em Proof.}
By hypothesis there is $\om \in \disc$ such that 
$L_{m}(\xi ) = L_{m}(\om )$. Hence, by Proposition \refth{compo}  
$L_{m} \circ L_{\xi} ( \lambda_{m, \xi} z)   =  
L_{m} \circ L_{\om} (\lambda_{m,\om} z)$ 
for every $z\in \disc$. The result follows because   
$L_m : \disc \rr P(m)$,  
$L_{\om}( \lambda_{m,\om}\, \_ ) : \disc \rr \disc$ and 
$L_{\xi}( \lambda_{m,\xi}\, \_ ) : \disc \rr P(\xi )$  are onto and one-to-one. 
%
\diaend

\begin{lemma}                \labeth{mul}
Let $\gamma \in {\Bbb C}$ with $|\gamma |=1$ and let $y\in \chu$. Then 
\begin{equation}    \laba{muly}     
 L_{\gamma y} (z) = \gamma L_{y} (\ov{\gamma}z)  . 
\end{equation} 
\end{lemma} 
{\em Proof.} 
Let $f\in \papa$ and $\{ z_{\alpha} \}$ be a net in $\disc$ tending to $y$. 
Thus $\gamma z_{\alpha}  \rr \gamma y$ and for every $z\in \disc$, 
\[
f(   L_{\gamma  y}(z)  )   =   \lim_{\alpha} f(  L_{\gamma z_{\alpha}}(z)  )  
					=   \lim_{\alpha} f(  \gamma L_{z_{\alpha}} (\ov{\gamma}z)  ) 
					=   f(  \gamma L_{y}(\ov{\gamma} z)  )  , 
\] 
as desired.      \diaend

\section{A characterization of Hoffman algebras} 
\equnew

\noindent 
{\sc Definition.} Let $m\in \calg$. The $m$-saturation of a set $E\subset \chu$ 
is defined as $L_{m}^{-1} ( L_{m} (E))$, and $E$ will be called $m$-saturated 
if it coincides with its $m$-saturation. We also write 
${\cal L}_{m}(y) \stackrel{\mbox{\scriptsize{def}}}{=}  
L_{m}^{-1} ( L_{m} (y))$ for 
the $m$-saturation of  $y\in \chu$.    \\

\noindent 
It is clear that ${\cal L}_{m}(0) \cap \disc = \{ 0 \}$. For $f\in \papa$ write 
\[   \zd (f) = \{ z\in \disc : f(z) =0 \}    \  \mbox{ and }\    
Z(f) = \{ \varphi \in \chu :  f(\varphi ) = 0 \} .  
\] 
It is well known that if $f$  is an interpolating Blaschke product then $Z(f)$ is the closure of $\zd (f)$. This immediately implies that if $m\in \calg$ and $y\in \chu$ 
are different points then there is an interpolating Blaschke product $f$ such that 
$f(m) = 0 \not = f(y)$. As a consequence we obtain that if  $m\in \calg$ then  
\[    {\cal L}_{m} (0) = \bigcap     \{  Z(b\circ L_{m})    :  
b  \mbox{ is an interpolating Blaschke product with } b(m) = 0  \} . 
\] 
Since $Z(b\circ L_{m})$ is the zero set of an interpolating Blaschke product then 
%
${\cal L}_{m}(0)$ is an intersection of closures of interpolating 
sequences, which in a sense is quite small.  Furthermore, 
\cite[Thm.$\,$1.4]{g-l-m} implies that $P(m)$ is a homeomorphic disk 
if and only if ${\cal L}_{m}(0) = \{  0  \}$. 

\begin{lemma}      \labeth{g-class}
Let $m\in \calg$. Then for $\om \in \disc$ we have 
\[  {\cal L}_{m}(\om ) =\{  L_{x}( \lambda_{m,x} \om ) : \,  x\in {\cal L}_{m}(0)  \} . 
\] 
\end{lemma} 
{\em Proof.} 
Since $L_{m}(x) = m$ for all $x\in {\cal L}_{m}(0)$ then by 
Proposition \refth{compo}
$L_{m} (L_{x} (\lambda_{m,x} z)) = L_{m} (z)$ on $\disc$. So, 
$L_{x}(\lambda_{m,x} \om )  \in {\cal L}_{m} (\om )$. 
Let $\xi \in \chu$ such that $L_{m}(\xi ) = L_{m}(\om )$. The `onto' part of 
Corollary \refth{o-o} implies that there is $x\in P(\xi )$ such that  $L_{m}(x) =m$. 
So,  $L_{m} \circ L_{x} (\lambda_{m,x} z) = L_{m}(z)$ for $z\in \disc$. 
In particular,  
$L_{m} ( L_{x} (\lambda_{m,x} \om ) ) = L_{m}(\om ) = L_{m}(\xi )$. 
That is, $L_{m}$ takes the same value on the points 
$L_{x}(\lambda_{m,x} \om )$ and $\xi$, which belong to $P(\xi )$. 
The `one-to-one' part of Corollary \refth{o-o} implies that 
$L_{x}(\lambda_{m,x} \om )=\xi$.    \diaend  


\begin{theo}    \labeth{big}
Let $m\in \calg \setminus \disc$ and let  $f\in \papa$ such that 
\[  f(   L_{x}(  \lambda_{m, x} z  )    ) = f(z)  
\ \mbox{ for all }\  x\in {\cal L}_{m}(0) \ \mbox{ and all }\  z\in \disc  . 
\] 
Then there is $F\in \papa$ such that $F\circ L_{m} = f$. 
\end{theo}

\noindent 
For $m\in \calg$ the theorem and Lemma \refth{g-class}  provide a description of  the 
algebra $\papa \circ L_{m}$  as 
the functions $f\in \papa$ such that  $f$  is constant on $L_{m}^{-1} (m')$  
for every $m' \in P(m)$.    
Only the sufficiency needs to be proved. 
We devote the next two sections to prove Theorem \refth{big}. 
For the sake of clarity      
it is convenient to rescue the hat for the Gelfand transform in the 
next corollary. 

\begin{coro}      \labeth{thm?} 
Let $m\in \calg \setminus \disc$ and $h   : P(m) \rr {\Bbb C}$ such that 
$h   \circ L_{m} = f    \in \papa$. Then the following conditions are equivalent.   
\begin{enumerate} 
\item[{\em (a)}] $h   $ is continuous on $P(m)$ with the topology induced by $\chu$,  
\item[{\em (b)}] $\hat{f   } \circ L_{x} (\lambda_{m,x} z)  =  f   (z)$   for every   
                $x\in {\cal L}_{m} (0)$   and    $z\in \disc$,    
\item[{\em (c)}] there exists $F\in \papa$ such that  
                $\hat{F} \circ L_{m} (z) = h   \circ L_{m} (z)$ for every $z\in \disc$, and  
\item[{\em (d)}] there exists $F\in \papa$ such that 
                $\hat{F} |_{P(m)} = h   $. 
\end{enumerate}  
\end{coro}  
{\em Proof.} 
We assume first that (a) holds. 
If (b) fails then there are 
$x\in {\cal L}_{m}(0)$ and $z_{0} \in \disc$ such that 
$\alpha = | \hat{f   } (L_{x} (\lambda_{m,x} z_{0}) ) -  f   (z_{0})|  > 0$. 
By the density of $\disc$ in $\chu$ the point  
$L_{x} (\lambda_{m,x} z_{0})$ is in the closure of the set 
$U = \{ z\in \disc : | \hat{f   } (L_{x} (\lambda_{m,x} z_{0}) )  -   f   (z)| < \alpha /2  \}$.         
Thus, by Lemma \refth{g-class}  
$  
L_{m}(z_{0})   =   L_{m} ( L_{x} (\lambda_{m,x} z_{0}) )   \in  
L_{m}(\ov{U})  \subset  \ov{L_{m}(U)}  . 
$ 
The continuity of $h   $ on $P(m)$ now implies that 
$h    ( L_{m}(z_{0}) )  \in \ov{  h   ( L_{m}(U) )  }$. 
But  this contradicts the fact  that  for  $z\in U$,  
\begin{eqnarray*} 
|h   (L_{m}(z_{0}))   -  h   (L_{m}(z))|    &  =  &    |f   (z_{0}) - f   (z)|     \\
&  \geq  &   
|f   (z_{0}) - \hat{f   } ( L_{x} (\lambda_{m,x} z_{0}) )|  -  
|\hat{f   } ( L_{x} (\lambda_{m,x} z_{0}) )  -   f   (z) |                \\
&  >  &             \alpha - \alpha /2  = \alpha  /2   .
\end{eqnarray*} 
The implication (b)$\Rightarrow$(c) is Theorem \refth{big}. 
Since $P(m) = L_{m} (\disc )$ then (d) is just a rephrasing of (c).  
Now suppose that (d) holds. 
Since $F\in \papa$ then $\hat{F}$ is continuous on $\chu$, and 
consequently $\hat{F}|_{ P(m) } = h   $ is continuous on $P(m)$. 
\diaend  

\section{Technical lemmas} 
\equnew 

The hyperbolic metric for $z,\om \in \disc$ is 
\[    \hh ( z, \om ) = \log \frac{1+\rho (z,\om )}{1-\rho (z,\om )}     . 
\] 
So, $\hh$ and $\rho$ are increasing functions of each other and 
$\hh (z, \om )$ tends to infinity if  and only if $\rho (z, \om )$ tends to $1$.  
We will use alternatively one metric or the other according to convenience. 
The hyperbolic ball of center $z\in \disc$ and radius $r>0$ will be denoted 
by  $\hypdi (z,r)$. 

%
The next two lemmas are easy consequences of  similar results in 
Hoffman's paper \cite[pp.$\,$82 and 86-88]{hof}. 
(or see \cite[pp.$\,$404-408]{gar}).  

\begin{lemma}       \labeth{a2}
Let $S$ be an interpolating sequence and let $m\in \ov{S}$. Then for any 
$0<\delta_{0}<1$ there is a subsequence $S'$ of $S$ such that 
$m\in \ov{S'}$ and $\delta (S') > \delta_{0}$.  
\end{lemma} 

%

\begin{lemma}    \labeth{dist} 
Let $b$ be an interpolating Blaschke product with $\delta (b) \geq \delta$ and let 
$\om \in \disc$. 
Then there is $0< c = c (\delta ) <1$      
such that  $c \rr 1$ as $\delta \rr1$, and    
\[   |b(\om )| \geq c    \rho (\om  , \zd (b)) . 
\] 
\end{lemma} 

\noindent 
Our next lemma is a trivial consequence of Lemma \refth{dist}; we state it 
for convenience.  

\begin{lemma}                   \labeth{hoffn}
Let $0< \sigma <1$. Then there are functions 
$0<  \delta (\sigma ) <1$    and   $s(\sigma ) >0$    
such that if $b$ is an interpolating Blaschke product with  
$\delta (b) > \delta(\sigma )$ 
then 
\[ |b(z)| \geq \sigma  \ \mbox{ for }\ 
z\not \in \bigcup \{ \hypdi  (z_{n}, s(\sigma )): \ z_{n} \in \zd (b) \} .
\]
\end{lemma}

\begin{lemma}                  \labeth{unif}
Let $\{ z_{n} \}$ be an interpolating sequence and 
$(z_{\alpha})$ be a subnet with  $z_{\alpha} \rr m \in \chu$. 
Suppose that  $F\in C(\chu )$. Then  
$F\circ L_{ z_{\alpha}} \rr  F\circ L_{m}$ 
uniformly on compact subsets of $\disc$.
\end{lemma} 
{\em Proof.} 
Since $z_{\alpha} \rr m$ in $\chu$ then $L_{z_{\alpha}} \rr  L_{m}$ in 
$\chu^{\Bbb D}$, and by the continuity of $F$ on $\chu$ then 
$F\circ L_{ z_{\alpha}} \rr  F\circ L_{m}   \stackrel{\mbox{\scriptsize{def}}}{=} f$ 
pointwise on $\disc$. 
We will see that  
 $F( L_{ z_{\alpha}}(z) )  \rr  f(z)$  uniformly on $|z| \leq r$ for any $0<r<1$. 
In fact, otherwise there is $\varepsilon >0$, a subnet  $(z_{\beta})$ of $(z_{\alpha})$ 
and points $\om_{\beta}$ with $|\om_{\beta}| \leq r$ such that 
\begin{equation}         \laba{far} 
| F(   L_{z_{\beta}} (\om_{\beta})  )  - f( \om_{\beta} ) |   > \varepsilon   
\ \ \mbox{ for all }\ \  \beta .
\end{equation} 
Taking a subnet of  $(z_{\beta})$ if necessary, we can also assume that 
$\om_{\beta} \rr \om$, with $|\om | \leq r$. 
Since $f$ is continuous at $\om$ and  
$F( L_{z_{\beta}} (\om ) ) \rr  f(\om )$, then 
there is $\beta_{0}$ such that for every $\beta \geq \beta_{0}$, 
\[   |f(\om ) - f(\om_{\beta})|  < \varepsilon /4  \ \ \mbox{ and }\ \  
|F( L_{z_{\beta}} (\om ) ) - f(\om )| < \varepsilon /4  . 
\] 
These inequalities together with (\ref{far})  give 
\[   
| F(   L_{z_{\beta}} (\om_{\beta})  ) -   F(   L_{z_{\beta}} (\om )  )|    > \varepsilon   /2   
\ \ \mbox{ for all }\ \  \beta \geq \beta_{0}  . 
\]
This will contradict the continuity of $F$ if we prove that 
$L_{z_{\beta}} (\om_{\beta})$ tends to $L_{m} (\om )$. 
Let $ (  L_{ z_{\gamma} } (\om_{\gamma})  )$ be an arbitrary convergent subnet 
of  $ (  L_{ z_{\beta} } (\om_{\beta})  )$, say to $y\in \chu$. 
Then by the lower semicontinuity of $\rho$, 
\[  \rho (y, L_{m}(\om ))   \leq     
\lim    \rho ( L_{ z_{\gamma} } (\om_{\gamma})  ,   L_{ z_{\gamma} } (\om ) )  = 
\lim    \rho (\om_{\gamma} , \om ) =0 , 
\]
meaning that $y= L_{m}(\om )$.  
So, every convergent subnet of   $ (  L_{ z_{\beta} } (\om_{\beta})  )$   
tends to $y$, 
and consequently the whole net tends to $y$.              \diaend \\ 

\noindent 
An immediate consequence of Lemma \refth{unif} is that if $F\in C(\chu )$ and 
$S$ is an interpolating sequence with $m\in \ov{S} \setminus S$, then for every 
$0<r<1$ and $\varepsilon >0$,  
\[   \left\{          z_{n} \in S : 
\sup_{|z| \leq r}   |F\circ L_{z_{n}} (z)  -  F\circ L_{m} (z)|  <  \varepsilon   \right\}  
\]
is a subsequence of $S$ having $m$ in its closure. 
 
\begin{lemma}                    \labeth{rhor}
Let $\xi , \om \in \disc$ and $m\in \chu \setminus \disc$.  
Then for any\/ $0<r<1$, 
\[  
\sup_{ |z|\leq r}  \rho ( L_{\xi} (\lambda_{m,\xi} z) , L_{\om} (\lambda_{m,\om}  z)  ) 
< \frac{30}{(1-r)^{2}} \rho (\xi ,\om ) .
\]
\end{lemma}
{\em Proof.} 
We can assume that $\pi (m) =1$. Also, 
since the desired inequality is obvious for 
$\rho (\xi ,\om ) > 1/30$, we can assume that $\rho (\xi ,\om ) < 1/2$. 
So, for $z\in \disc$ with $|z|\leq r$, 
\begin{eqnarray*} 
\rho ( L_{\xi} (\lambda_{m,\xi} z) , L_{\om} (\lambda_{m,\om} z)  )  
&  \leq   &  
\rho ( L_{\xi} (\lambda_{m,\xi} z) , L_{\xi} (\lambda_{m,\om} z)  )     +     
\rho ( L_{\xi} (\lambda_{m,\om} z) , L_{\om} (\lambda_{m,\om} z)  )  \\ 
&   =      &     \varrho_{1} + \varrho_{2} .
\end{eqnarray*} 
Since $L_{\om} : \disc \rr \disc$ is an onto isometry with respect to $\rho$, 
then there exists $v\in \disc$ 
such that $\xi = L_{\om}(v)$ and $|v| = \rho (\xi , \om )$. 
The elementary formula (for $\pi (m) =1$) 
\[                                        
\lambda(m, L_{\om}(v))      =     
\frac{1+ L_{\om}(v)}{1+ \ov{L_{\om}(v)}}   =   
\left( \frac{1+\om \ov{v}}{1+\ov{\om}v} \right)  
\left( \frac{ 1+v \ov{ \lambda (m,\om )} }{ 1+ \ov{v}  \lambda (m,\om ) } \right)  
\lambda (m,\om )  
\]
yields   
\begin{eqnarray}            \label{t38} 
|\lambda (m, L_{\om}(v))  -  \lambda (m,\om )|   &  =  & 
\frac{  | v(  \ov{\lambda}- \ov{\om} )  +  \ov{v}( \om -\lambda )+     
|v|^{2}( \om \ov{\lambda} - \ov{\om}\lambda )|  }{ 
 | (1+\ov{\om}v)  (1+\ov{v} \lambda )  |  }  \nonumber     \\ 
 &  <  & \frac{6|v|}{ (1-|v|)^{2} }  < 24 |v|  ,
\end{eqnarray} 
where $\lambda = \lambda (m,\om )$, and the last inequality holds because $|v| < 1/2$. 
Thus, (\ref{t38}) and the isometric property of  $L_{\xi}$  give 
\begin{equation}   \laba{r1}
\varrho_{1} = \rho (  \lambda_{m,\xi}z , \lambda_{m,\om}z )  \leq 
\frac{ |\lambda (m,\xi ) - \lambda (m,\om )| }{ 1-r^{2} } < 
\frac{24 |v|}{ 1-r^{2} } < \frac{24 |v|}{ (1-r)^{2} }    . 
\end{equation} 
Put $z' = \lambda_{m,\om}z$. 
Then by (\ref{pong})
\begin{eqnarray*}
\varrho_{2} 
&   =   &    \rho (  L_{ L_{\om}(v) } (z') , L_{\om}(z')  )   
\     =   \      \rho (  L_{\om} \circ L_{v} (  \lambda_{\om ,v}z' )  ,  L_{\om}(z')  ) \\[1.5mm]  
&   =   &    \rho (  L_{v} (  \lambda_{\om ,v}z' )  ,  z'  ) 
\   \leq  \     \frac{| L_{v} (\lambda_{\om ,v} z' )  -  z' |}{ 1-|z'| }   \\[1mm]  
&   =   &   \frac{ |z' (\lambda_{\om ,v} -1) + 
v -\ov{v}  (z')^{2}  \lambda_{\om ,v} | }{ (1-|z'|)  \,  |1+\ov{v} \lambda_{\om , v} z'|  }     
\   \leq  \     \frac{ | \lambda_{\om ,v} -1| + 2|v| }{(1-r)^{2}} .
\end{eqnarray*} 
This inequality together with 
\[  
| \lambda_{\om ,v} -1| = \frac{ |\ov{\om} v - \om \ov{v}| }{|1+\om \ov{v}|} 
< \frac{2|v|}{ 1-|v| }  < 4|v| 
\] 
yields  $\varrho_{2} < 6|v| (1-r)^{-2}$. So, adding this estimate to (\ref{r1}) we 
obtain 
\[   \varrho_{1} + \varrho_{2} <  (24 + 6) |v|  (1-r)^{-2}   =
30 \rho (\xi , \om )  (1-r)^{-2}, \] 
as promised.  
\diaend

\begin{lemma}                    \labeth{sat}
Let $m\in \calg$ and  $E\subset \chu$ be a closed $m$-saturated set. 
If $V$ is an open neighborhood of $E$ then there exists an open  
$m$-saturated set\/ $W$ such that  $E\subset W\subset V$. 
\end{lemma}
{\em Proof.}
Since $\chu \setminus V$ is closed then so is $L_{m}(\chu \setminus V)$. 
Since $E$ is $m$ saturated then the closed $m$-saturated set 
$F\stackrel{\mbox{\scriptsize{def}}}{=}  
L_{m}^{-1} ( L_{m}(\chu \setminus V) )$ does not meet $E$. 
Additionally, $F\supset \chu \setminus V$ and then the open set 
$W \stackrel{\mbox{\scriptsize{def}}}{=} \chu \setminus F$ satisfies the lemma. 
\diaend  

\begin{lemma}                                 \labeth{sat2}
Let $S$ be an interpolating sequence and $m\in \ov{S}$. 
If\/ $W\subset \chu$ is an open $m$-saturated neighborhood of  ${\cal L}_{m}(0)$ 
then there is a subsequence $S_{0} \subset S$ such that 
$m\in \ov{S}_{0} \cap \ov{P(m)} \subset L_{m}(W)$.
\end{lemma}
{\em Proof.}
The hypothesis on $W$ implies that $L_{m}(W)$ is open in $\ov{P(m)}$ and 
$m\in L_{m}(W)$. That is, the set 
$E =    
\ov{P(m)}\setminus L_{m}(W)$ is 
closed and $m\not \in E$. 
By compactness then there is an open set $U\subset \chu$ such that $m\in U$ and 
$\ov{U}\cap E= \emptyset$. 
Defining $S_{0}= S\cap U$ we have that $m\in \ov{S}_{0} \subset \ov{U}$. 
Hence, $\ov{S}_{0} \cap E = \emptyset$ and then 
$\ov{S_{0}} \cap \ov{P(m)}  \subset   \ov{P(m)} \setminus E = L_{m}(W)$. 
\diaend  

\begin{lemma}                           \labeth{seven} 
Let $m\in \calg$ and $f\in \papa$ such that 
$f\circ L_{x} (\lambda_{m,x} z) = f(z)$ for all $x\in {\cal L}_{m}(0)$  and  $z\in \disc$. 
For $\varepsilon >0$ and $0<r<1$ consider the set 
\[    
U= \{  \om \in \disc :  | f\circ L_{\om} (\lambda_{m,\om} z) - f(z)|  <  
\varepsilon     \ \mbox{ for }\     |z| \leq r   \}  .
\]    
Then $\ov{U}$ is a neighborhood of   ${\cal L}_{m}(0)$.  
\end{lemma}
{\em Proof.} 
If the lemma fails then there is $x\in {\cal L}_{m}(0)$ in the closure of  
$V= \chu \setminus \ov{U}$. 
Since $V$ is open and $\disc$ is dense in $\chu$, a simple topological 
argument shows that 
$\ov{V} = \ov{ V \cap \disc }$. Therefore $x\in  \ov{ V \cap \disc }$, where 
\[    V \cap \disc     \subset   
 \{  \om \in \disc :  | f\circ L_{\om} (\lambda_{m,\om} z_{\om}) - f(z_{\om})|  \geq     
\varepsilon     \ \mbox{ for some }\   
z_{\om} \ \mbox{ with }\    |z_{\om}| \leq r   \}  .
\]
Let $(\om_{\alpha})$ be a net in $V\cap \disc$ that tends to $x$, and write  
$z_{\alpha} \stackrel{\mbox{\scriptsize{def}}}{=} z_{\om_{\alpha}}$. 
By taking a suitable subnet we can also assume that $z_{\alpha} \rr z_{0}$, 
where $|z_{0}| \leq r$. 
Thus, 
\begin{equation}        \laba{ripi} 
| f\circ L_{\om_{\alpha}} (\lambda_{m,\om_{\alpha}} z_{\alpha}) - f(z_{\alpha})| 
\geq   \varepsilon    \ \mbox{ for every }\     \alpha      . 
\end{equation} 
We can assume $\| f \|_{\infty} \leq 1$. 
By the Schwarz-Pick inequality \cite[p.$\,$2]{gar},     
\begin{eqnarray*} 
\rho ( f\circ L_{\om_{\alpha}} (\lambda_{m,\om_{\alpha}} z_{\alpha}) , 
         f\circ L_{\om_{\alpha}} (\lambda_{m,x} z_{0})  ) 
&  \leq  &   
\rho (  \lambda_{m,\om_{\alpha}} z_{\alpha}  ,   \lambda_{m,x} z_{0}  )  \\
&  \leq  & 
\rho (  \lambda_{m,\om_{\alpha}} z_{\alpha}  ,   \lambda_{m,x} z_{\alpha}  )  + 
\rho (  \lambda_{m,x}  z_{\alpha}  ,   \lambda_{m,x} z_{0}  )    \\
&  \leq  & 
\frac{1}{1-|z_{\alpha}|^{2}}  | \lambda_{m,\om_{\alpha}}   -    \lambda_{m,x} |   
+  \rho (z_{\alpha} , z_{0})  , 
\end{eqnarray*}
which tends to zero.  
Since $L_{\om_{\alpha}} \rr  L_{x}$ then the last inequality gives 
\[ 
\lim_{\alpha} f\circ L_{\om_{\alpha}} ( \lambda_{m,\om_{\alpha}} z_{\alpha} )  = 
\lim_{\alpha} f\circ L_{\om_{\alpha}} ( \lambda_{m,x} z_{0} )  =
                     f\circ L_{x} ( \lambda_{m,x} z_{0} )  =   f(z_{0})  ,  
\]
which contradicts (\ref{ripi}).    \diaend  \\   

\noindent 
The next lemma is in \cite[Lemma 2.1]{su4}.    

\begin{lemma}                                   \labeth{useful} 
Let $u$ be an inner function and $0<\beta <1$. Put 
$V = \{ z\in \disc : |u(z)| < \beta \}$ and suppose that $f\in \papa (V )$.   
Then there are  $0< \gamma = \gamma (\beta ) <\beta$, 
$C= C(\beta )>0$ and $F\in \papa$ such that  
%
\begin{enumerate} 
\item[{\em (i)}] $\| F \|_{\infty} \leq C \| f  \|_{ \papa (V ) } $, 
and  
%
\item[{\em (ii)}] $|F(z) - f(z)|         \leq           A  \| f  \|_{ \papa (V ) } \/   |u(z)| \,$  when  
$\, |u(z)| < \gamma$,  where $A= \gamma^{-1}(C+1)$.  
\end{enumerate} 
\end{lemma} 

\section{Proof of Theorem \refth{big}}  
\equnew

Given $m\in \calg \setminus \disc$ and $f\in \papa$ that satisfy the hypotheses 
of the theorem, we are going to construct a function $F\in \papa$ 
such that $F\circ L_m = f$. 
We can assume without loss of generality 
that $\pi (m) =1$ and $\| f \| =1$. 
Let $\{ \sigmak \} \subset (0,1)$ be a sequence satisfying  
$\prod_{k\geq 1} \sigmak  > 0$ and let $s(\sigmak )$ be the associated 
parameters given by Lemma \refth{hoffn}. 
Take $s_{k} > s(\sigmak )$ tending 
{\em increasingly} to $\infty$,  
and put  $r_{k}  \stackrel{\mbox{\scriptsize{def}}}{=} 
4\sspace  (2^{k} s_{1} + 2^{k-1} s_{2} + \cdots + 2 s_{k})$. 

Given an arbitrary interpolating sequence $S$ such that  $m\in \ov{S}$, and 
$\{ \vark \}  \subset  (0,1)$ a decreasing sequence that tends to 0, we will 
construct a decreasing chain of subsequences $S_{k} = \{ z_{k,n} : n\geq 1 \}$, 
$S \supset S_{1} \supset S_{2} \supset  \cdots$,
such that for every $k\geq 1$,                   
\mbox{ } \vspace{1.5mm}  \mbox{ }\\
\noindent 
 (1)  $m\in \ov{S}_{k}$,  \\
 (2)  $\hh  (      \{ |z-1|   >  \vark \}  , S_{k}    ) > r_{k}$, \\ 
 (3)  $\sum_{n\geq 1} (1-|z_{k,n}|)  < 2^{-k}$,    \\ 
 (4)  $\hh ( z_{k, n_{1}}, z_{k,n_{2}}) > r_{k} \ \mbox{ for }\ 
      n_{1} \not = n_{2}$, \\ 
 (5)  if $b_{k}$ is the interpolating Blaschke product with zero sequence $S_{k}$, then 
      \[  |b_{k}(z)| \geq \sigmak   \ \ \mbox{ for }\ \   z\not \in 
      \bigcup_{n\geq 1} \hypdi (z_{k,n}, s_{k})  ,  
      \]  
 (6)  if   $l\geq k$  and  $\hh ( z_{l,p}, z_{k,n}) < r_{l}$  then 
      $\rho (  L_{-z_{l,p}}(z_{k,n})  ,  T_{k} ) < \vark$,  
      where   
      $\, T_{k} =  L_{m}^{-1} (   \ov{S}_{k} \cap \ov{P(m)}   )$,  and \\ 
 (7)  $| f\circ L_{ \om  } (  \lambda_{m,\om} z  ) - f(z) | \leq  \vark$    
       \mbox{ when } 
      $\hh (z,0) \leq s_{k}$    and    $\om \in T_{k} \cap \disc$.  \\ 

\noindent  
{\bf The first construction.} 
The argument will be inductive. 
By Lemmas \refth{seven} and \refth{sat},  
for every $k\geq 1$ there is an open $m$-saturated 
neighborhood of  ${\cal L}_{m}(0)$,   $W_{k} \subset \chu$,   
such that $W_{k+1} \subset W_{k}$ and for all $\om \in W_{k}$, 
\begin{equation}      \laba{qu1}
  | f\circ L_{ \om  } (  \lambda_{m,\om}z  ) - f(z) | \leq \vark    
       \ \mbox{ if  }\       \hh  (z,0) \leq s_{k}. 
\end{equation} 
Step 1. By Lemma \refth{sat2} there is $S'_{1} \subset S$ such that (1) holds and 
$\ov{S'_{1}} \cap \ov{P(m)}     \subset     L_{m}(W_{1})$. 
Since $W_{1}$ is $m$-saturated then  
$T'_{1}\stackrel{\mbox{\scriptsize{def}}}{=}  
L_{m}^{-1}(   \ov{S'_{1}} \cap \ov{P(m)}   )  \subset 
L_{m}^{-1} (L_{m} (W_{1})) = W_{1}$. 
Then (\ref{qu1}) tells us that $S'_{1}$ satisfies (7), and then 
so does any subsequence of  $S'_{1}$ that contains $m$ in its closure. 
By Lemmas \refth{a2} and \refth{hoffn} 
we can assume that $\delta (S'_{1})$ is so close to 1 
that (4) and (5) hold. 
Furthermore, since $\pi (m) =1$ we can easily achieve 
conditions (2) and (3) by taking as $S_{1}$ the subsequence of $S'_{1}$ whose 
elements  are contained in a sufficiently small Euclidean ball centered at $1$. 

Condition (6) only makes sense for $k=l=1$.  
If   $z_{1,n} , z_{1,p} \in S_{1}$   are such that   $\hh (z_{1,n} , z_{1,p}) < r_{1}$  
then (4) implies that   $z_{1,n} = z_{1,p}$. 
Therefore $L_{-z_{1,p}}  (z_{1,n}) = 0 \in  T_{1} \cap \disc$, 
because $L_{m}(0) = m \in \ov{S}_{1} \cap \ov{P(m)}$.  \\ 

\noindent
Step $l$. Let $l\geq 2$ and suppose that we already have 
$S\supset S_{1} \supset  \ldots \supset S_{l-1}$ satisfying 
$(1) \ldots (7)$. 
By Lemma \refth{sat2} there exists $S'_{l} \subset S_{l-1}$ 
such that (1) holds and $\ov{S'_{l}} \cap \ov{P(m)} \subset L_{m} (W_{l})$. 
Since $W_{l}$ is $m$-saturated then 
$ T'_{l} \stackrel{\mbox{\scriptsize{def}}}{=}  
L_{m}^{-1}( \ov{S'_{l}} \cap \ov{P(m)} )  \subset  
 L_{m}^{-1} (L_{m} (W_{l}))  =   W_{l}$.  
By (\ref{qu1}) the sequence $S'_{l}$ satisfies (7), and the same holds 
for any subsequence of $S'_{l}$ having $m$ in its closure. 
As in the case $l=1$, by Lemma \refth{a2} and \refth{hoffn} 
we can assume that $S'_{l}$ 
satisfies (4) and (5), and by taking the points of $S'_{l}$ that are close enough 
to $1$ we can also assume that $S'_{l}$ satisfies (2) and (3). 
Clearly, any subsequence $S_{l}$ of  $S'_{l}$ such that $m\in \ov{S_{l}}$  will 
satisfy all the above properties. Therefore we will be done if we can pick the 
sequence $S_l$ so that it also satisfies (6).

Let $k\leq l-1$ and let $\eta_{k}>0$ to be chosen later. 
By \cite[Lemma 1.8]{g-l-m}  we have $b_{k} \circ L_{m} = B_{k} g_{k}$, 
where $g_{k} \in (\papa )^{-1}$ and  
$B_{k}$ is an interpolating Blaschke product with zero sequence  
$Z_{\disc}(B_{k})$ contained in $T_{k}$. 
The inclusion holds because if $B_k(z_0)=0$ then $b_k(z_0)=0$, 
and consequently 
$L_m(z_0) \in P(m) \cap Z(b_k) = P(m) \cap \ov{S}_k$. 

Since $m\in \ov{S'_{l}}$ then by 
the remark following Lemma \refth{unif} 
there is a  subsequence  $\Lambda_{k} \subset S'_{l}$ such that  
$m\in \ov{\Lambda}_{k}$  and 
\begin{equation}   \laba{qu2} 
| b_{k} \circ L_{\nu} (z)  - B_{k}(z) g_{k}(z) | < \eta_{k} 
\ \mbox{ for }\           \nu \in \Lambda_{k}    \ \mbox{ and }\  \hh (z,0) < r_{l}. 
\end{equation} 
If  $\nu \in \Lambda_k$ and $z_{k,n}\in S_{k}$     satisfy     
$\hh (z_{k,n} , \nu ) < r_{l}$     then  
$\hh (L_{-\nu}(z_{k,n}) , 0 ) < r_{l}$.  
Applying (\ref{qu2})  to $z = L_{-\nu}(z_{k,n})$ we get  
\[   |B_{k} (  L_{-\nu}(z_{k,n})  )|   \    |g_{k} (L_{-\nu}(z_{k,n}))|    =     
     | b_{k} \circ L_{\nu} (L_{-\nu}(z_{k,n})) -  (B_{k} g_{k})  (L_{-\nu}(z_{k,n}))|  
< \eta_{k} .
\]
We are using here that $L_{-\nu} = L_{\nu}^{-1}$ and  $b_{k}(z_{k,n}) =0$. 
Then  $|B_{k} (L_{-\nu}(z_{k,n}))| <  \eta_{k}  \| g_{k}^{-1} \|_{\infty}$. 
Since $B_{k}$ is interpolating, Lemma \refth{dist} implies that for 
small values of $\eta_{k}$ the point $L_{-\nu}(z_{k,n})$ must be close to 
the zero sequence of $B_{k}$ in the $\rho$-metric. 
That is, choosing $\eta_{k}$ small enough we obtain 
\begin{equation}    \laba{qu3} 
  \rho (      L_{-\nu}(z_{k,n}) ,  T_{k}  )   \leq  
    \rho (      L_{-\nu}(z_{k,n}) ,  Z_{\disc} (B_{k})    )   < \vark 
\end{equation} 
for every $\nu \in \Lambda_{k}$
such that $\hh (\nu , z_{k,n}) < r_{l}$ for some $z_{k,n}$.  
Doing    
this process for $k=1, \ldots , l-1$ we obtain the respective subsequences $\Lambda_{k}\subset S'_{l}$ 
satisfying (\ref{qu3}), and such that $m\in \ov{\Lambda}_{k}$ 
for $k=1, \ldots , l-1$.  
Since disjoint subsequences of an interpolating sequence have disjoint closures 
then $m$ is in the closure of  
\[  S_{l} \stackrel{\mbox{\scriptsize{def}}}{=} 
 \bigcap_{1\leq k\leq l-1} \Lambda_{k} ,  
\] 
and by (\ref{qu3}) $S_{l}$ satisfies (6) for $k= 1,\ldots , l-1$. 
Finally, the same argument used in Step 1 
shows that $S_{l}$ satisfies (6)  also for $k=l$. \\

\noindent 
{\bf The second construction.} 
Observe that condition (4) implies that for a fixed value of $k$, 
\[  \hypdi (z_{k, n_{1}}, s_{k}) \cap \hypdi (z_{k, n_{2}} , s_{k}) = \emptyset  
\ \ \mbox{ if }\ \  n_{1} \not = n_{2} .
\]
Now we define recursively some sets made of unions of the balls 
$\hypdi ( z_{k,n}, s_{k})$,   
which we call `swarms'. 
For $n\geq 1$ the swarm of height 1 and center $z_{1,n}$ is defined as 
\[  E_{1,n}=  \hypdi ( z_{1,n}, s_{1}) .
\]
Once we have the swarms of height $j= 1, \ldots , k-1$, 
we define the swarm of height $k$ and center $z_{k,n}$ 
(for $n\geq 1$) as 
\[ 
E_{k,n} =  \hypdi (z_{k,n} , s_{k})  \bigcup \{ E_{j,p} : \ j\leq k-1, \ p\geq 1 
\ \mbox{and}\  E_{j,p} \cap \hypdi (z_{k,n}, s_{k}) \not = \emptyset  \}     . 
\] 
We write $\diamh E = \sup \{ \hh (x,y): \, x, y\in E \}$ for the 
hyperbolic diameter of a set $E\subset \disc$. 
The next three properties will follow by induction. \\
\mbox{}\hspace{1.4mm} (I)    $\diamh E_{k,n}   \leq   2^{k} s_{1} + 2^{k-1} s_{2} 
+ \cdots + 2s_{k}$, \\
\mbox{}\hspace{0.11mm}  (II)   $E_{k,n_{1}}  \cap   E_{k,n_{2}}  = \emptyset \ $ 
 if  $\ n_{1} \not = n_{2}$, and \\ 
(III)  each swarm of height $j\leq k-1$ meets (and then it is contained in) at most 
one swarm of  height $k$. 

\noindent 
{\em Proof of {\em (I)}}. This is trivial for $k=1$. 
By the definition of swarms and inductive hypothesis,    
\begin{eqnarray*} 
\lefteqn{  \diamh E_{k,n}  \leq   }  \nonumber \\  
                                    & \leq   &   \diamh \hypdi (z_{k,n}, s_{k})      +    
   2\max 
\{ \diamh E_{j,p} : \, j\leq k-1, \ E_{j,p} \cap \hypdi (z_{k,n}, s_{k}) \not = \emptyset \}  \\
			  &  \leq   &   2s_{k} + 2\sspace ( 2^{k-1} s_{1} 
				       + 2^{k-2} s_{2} + \cdots + 2s_{k-1} )   
			    =            2^{k} s_{1} + \cdots + 2s_{k} .   
\end{eqnarray*} 

\noindent  
{\em Proof of {\em (II)} and {\em (III)}}.  
Suppose that $E_{j,p}$ (with $j\leq k$)  meets $E_{k,n_{1}}$ and 
$E_{k,n_{2}}$, where $n_{1} \not = n_{2}$. Then by (I) 
\begin{eqnarray*}  
\hh  (z_{k,n_{1}}, z_{k,n_{2}})  
 &  \leq   &  \diamh E_{k,n_{1}}  + \diamh E_{j,p} + \diamh E_{k,n_{2}}   \\
 &  \leq   &  3\sspace (2^{k} s_{1} + 2^{k-1} s_{2} + \cdots + 2s_{k} )  < r_{k} , 
\end{eqnarray*} 
which contradicts condition (4). When $j=k$ this proves (II), and for $j< k$ this proves 
(III) except for the statement between brackets.  
So, suppose that $E_{j,p}$  meets  $E_{k,n}$, where $j\leq k-1$. 
If  $E_{j,p}$ meets $\hypdi (z_{k,n} , s_{k})$ then 
$E_{j,p} \subset E_{k,n}$ by definition. Otherwise there is some 
swarm $E$ of height at most $k-1$ such that 
\[  E\cap \hypdi  (z_{k,n}, s_{k}) \not = \emptyset \ \ \mbox{ and }\ \ 
    E\cap E_{j,p} \not = \emptyset  . 
\]
If height $E\geq j =$ height $E_{j,p}$ then by inductive hypothesis (II) for the equality 
and (III) for the strict inequality, we have $E_{j,p} \subset E$. Hence, $E_{j,p}$  
is contained in $E_{k,n}$. 
Similarly, if height $E< j$ then  inductive  hypothesis (III) implies that 
$E\subset E_{j,p}$. Therefore 
\[ E_{j,p} \cap \hypdi ( z_{k,n} , s_{k}) \supset E \cap \hypdi ( z_{k,n} , s_{k}) 
\not = \emptyset , 
\]
and then $E_{j,p} \subset E_{k,n}$  by definition. 

Some remarks are in order. 
Condition (II) says that two 
swarms of the same height are either the same (with the same $n$) or they are disjoint, 
and condition (III) says that if two different swarms have non-void intersection, 
then the one of smaller height is contained into the other.  
Also, observe that by (II), $E_{k,n} \cap S_{k} = \{ z_{k,n} \}$. 

We will see that if  $l\geq k$ then 
\begin{equation}            \laba{s-b} 
\{  |z-1| > \varepsilon_{k}  \} \cap E_{l,p} = \emptyset  
\ \mbox{ for all }\   p\geq 1  . 
\end{equation} 
In fact, suppose that for some $l\geq k$ and $p\geq 1$ there is  
$\om \in E_{l,p}$ with $|\om  -1|> \vark$. 
Then, since $\{ \varepsilon_{j} \}$ is a decreasing sequence, 
(2) and  (I)  yield 
\[    
   h (\om , z_{l,p})   \geq    \hh  ( \{ |z-1| > \vark  \} , S_{l})  \\
		         \geq     \hh  ( \{ |z-1| > \varepsilon_{l} \} ,  S_{l}) \\
		           >      r_{l}   \geq     4 \sspace \diamh \, E_{l,p}  , 
\]
which is not possible. 
By (\ref{s-b}) every strictly increasing chain of swarms 
$E^{ (1) } \subset E^{ (2) } \subset \ldots$ is finite. Because if $\om \in E^{ (1) }$ then  
there is $k$ such that $|\om -1| > \vark$  (since $\varepsilon_{j}$ tends to 0), and 
therefore $E^{ (1) }$ cannot lie in any swarm of height  $\geq k$. 
Roughly speaking, we could say that there is no swarm of infinite height. 
Consequently, every swarm is contained in a unique maximal swarm, and 
\[  \Omega  \stackrel{\mbox{\scriptsize{def}}}{=}  
\bigcup_{k, n \geq 1}  \hypdi  (z_{k,n} , s_{k}) = 
\bigcup \{ E_{l,p}  \mbox{ maximal}   \}  . 
\]

\noindent
{\bf Choosing $\varepsilon_{k}$}.  
Define a function $g\in \papa (\Omega )$ by 
$g(\om ) = f\circ L_{z_{l,p}}^{-1} (\om )$ for $\om \in E_{l,p}$, with 
$E_{l,p}$ a maximal swarm. 
The only requirements that we have imposed so far to the sequence  
$\{ \varepsilon_k \}$ are that it is contained in $(0,1)$ and decreases to zero. 
We claim that there is a choice of the sequence $\{  \varepsilon_{k}  \}$ so that 
\begin{equation}      \laba{ely} 
    \lim_{k}   g( L_{z_{k,n}}(z))  = f(z) , 
\end{equation} 
where the limit is uniform on $n$ and on compact subsets of $\disc$. 
Fix $0<r<1$ and let $z\in \disc$ with $|z| \leq r <1$.  
Since $\lim_{k} s_{k} = \infty$ then for 
$k$ big enough   
we have 
\begin{equation}      \laba{ine1} 
|z|\leq  r  < (e^{s_{k}} - 1) / (e^{s_{k}}+1)  . 
\end{equation} 
This means that  $z\in \hypdi (0, s_{k})$.      
The point $z_{k,n}$ is in some maximal swarm $E_{l,p}$ with $l\geq k$. 
Hence by (\ref{ine1})  
\[    
L_{z_{k,n}} (z) \in \hypdi  (z_{k,n}, s_{k}) 
\subset E_{l,p} \subset \Om   
\]    
Since $g$ is defined on $\Om$ then $g( L_{z_{k,n}} (z)  )$ makes sense, and 
\begin{equation}     \laba{qu6}
g(L_{z_{k,n}} (z) ) =                                   
f\circ L_{z_{l,p}}^{-1}( L_{z_{k,n}}(z)  ) = 
f\circ L_{   L_{-z_{l,p}}(z_{k,n})   }  (  \ov{ \lambda (-z_{l,p} , z_{k,n}) } z   )  , 
\end{equation} 
where the last equality comes  from 
the identity $ L_{z_{l,p}}^{-1}  = L_{-z_{l,p}}$ and (\ref{pong}).       
A simple calculation shows that 
$ \ov{  \lambda (-z_{l,p} , z_{k,n})   }=  
      \lambda (  z_{k,n} , L_{-z_{l,p}}(z_{k,n})    )  . 
$  
So, if   $\xi \stackrel{\mbox{\scriptsize{def}}}{=} L_{-z_{l,p}}(z_{k,n})$ we 
can  rewrite (\ref{qu6}) as 
\begin{equation}     \laba{qu7} 
g(L_{z_{k,n}} (z) ) =                                   
f\circ L_{\xi}   (   \lambda_{z_{k,n} , \xi}   z   )  .   
\end{equation}
%
Since $z_{k,n} , z_{l,p} \in E_{l,p}$ then by (I),                    
$h (z_{k,n} , z_{l,p}) \leq  \diamh E_{l,p} < r_{l}$.  
Thus, (6) implies that there is 
$\om \in T_{k}$ such that  $\rho ( \xi , \om ) < \vark$. 
Since $h(z,0) \leq s_{k}$ (by (\ref{ine1})) and $\| f \|_{\infty} =1$,  
then successive applications of (\ref{qu7}),  (7) and the Schwarz-Pick inequality  
yield   
\begin{eqnarray}            \laba{qu8}
  |g( L_{z_{k,n}} (z) ) - f(z) |    
 &  \leq &  
 |f\circ L_{\xi}(   \lambda_{z_{k,n} , \xi}   z   ) - 
		   f\circ L_{\om}  ( \lambda_{m,\om} z )|    \nonumber  \\ 
                                      &   +   & 
  		  | f\circ L_{\om}  ( \lambda_{m,\om} z )  -  f(z)|   \nonumber \\  
			    &  \leq  &   
     2  \rho  (    f\circ L_{\xi}(   \lambda_{z_{k,n} , \xi}   z   ) ,  
		   f\circ L_{\om}  ( \lambda_{m,\om} z )   )    +   \vark   \nonumber \\
			   &  \leq  &   
     2  \rho  ( L_{\xi}(   \lambda_{z_{k,n} , \xi}  z   ) ,  L_{\om}( \lambda_{m,\om} z )   )   
			     +   \vark      , 
\end{eqnarray} 
where 
\begin{eqnarray*}
\rho  ( L_{\xi}( \lambda_{z_{k,n} , \xi}  z ) ,  L_{\om}( \lambda_{m,\om} z ) )
&  \leq   &
\rho  ( L_{\xi}( \lambda_{z_{k,n} , \xi}  z ) ,  L_{\xi}( \lambda_{m ,\xi} z) )  
  +   
\rho  (L_{\xi}( \lambda_{m ,\xi}  z)  ,  L_{\om}( \lambda_{m,\om} z )   )  \\
&   =     &  
\varrho_{1} + \varrho_{2}  . 
\end{eqnarray*} 
Using the isometric property of  $L_{\xi}$, a straightforward calculation shows that 
\begin{eqnarray}       \laba{qu11}  
\varrho_{1}  
&   =    & 
\rho  ( \lambda_{z_{k,n} , \xi}  z , \lambda_{m ,\xi}  z  ) 
\   \leq   \ 
\frac{ |\lambda_{z_{k,n} , \xi}  -   \lambda_{m ,\xi}| }{1-|z|^{2}}  \nonumber  \\
&   =   & 
\left|  \frac{   (1-z_{k,n}) (1+\xi ) \ov{\xi} -  (1-\ov{z}_{k,n}) (1+\ov{\xi} ) \xi    }{  
         (1-|z|^{2}) (1+\ov{\xi} z_{k,n}) (1+\ov{\xi})  }  \right|   \nonumber  \\
&    \leq   & 
\frac{2}{(1-r^{2})}   \frac{|1-z_{k,n}|}{ (1-|\xi |) }  . 
\end{eqnarray}
By (2) and since $\varj \rr 0$, there is $j> 1$ such that 
\begin{equation}    \laba{qu12}
\varj < |1-z_{k,n}| \leq \varepsilon_{j-1}  , 
\end{equation} 
and since $z_{k,n} \in E_{l,p}$ then (\ref{s-b}) 
says that $l = \mbox{height} \, E_{l,p} < j$. 
By (I) then $h (z_{k,n}, z_{l,p})  \leq  \diamh E_{l,p} <   r_{l} <  r_{j}$. 
Thus,  
\[
|\xi | = |L_{-z_{l,p}}(z_{k,n}) | = \rho ( L_{-z_{l,p}}(z_{k,n}) , 0 ) = 
\rho (z_{k,n}, z_{l,p}) < 
\frac{ e^{r_{j}} - 1}{ e^{r_{j}} + 1}  , 
\]
and consequently $(1-|\xi |) \geq  e^{ - r_{j} }$.   
Choosing 
$\varepsilon_{q-1}   =        (2 \, q \, e^{r_{q}})^{-1}$   for all   $q>1$,  
we can insert the last inequality and (\ref{qu12}) in (\ref{qu11}), thus obtaining  
\[          
\varrho_{1} 
\leq       \frac{ 2 e^{r_{j}} \,   \varepsilon_{j-1} }{  1-r^{2} } 
=     \frac{ j^{-1} }{  (1-r)^{2} }
<      \frac{ k^{-1} }{  (1-r)^{2} } .  
\]         
The last inequality holds because $j >  l  \geq  k$.  
On the other hand, since $|z| \leq r$ and $\rho (\xi , \om )< \vark$ 
then Lemma \refth{rhor} says that $\varrho_{2} < (1-r)^{-2} 30 \vark$.   
Putting all this together in (\ref{qu8})
we obtain that whenever $|z|\leq r$ and $k$ is large enough so that 
(\ref{ine1}) holds, then  
\[ 
  |g(L_{z_{k,n}} (z)) - f(z) |   \leq  
 2(\varrho_{1} + \varrho_{2}) + \vark 
<  \frac{C}{(1-r)^{2}} \,  \frac{1}{k}      
\] 
for some absolute constant $C>0$. 
Since this estimate is independent of $n$ then (\ref{ely}) follows. \\

\noindent 
{\bf The construction of $F$}.  
We recall that $b_{k}$ is a Blaschke product with zero sequence $S_{k}$. 
Since $S$ is an interpolating sequence, then  
$a= \inf \{ |b_{k}(z)| : z\in S\setminus S_{k} \} > 0$, 
and since $m\in \ov{S}_{k}$ then  $\{ x\in \chu : |b_{k}(x)| < a/2 \}$ 
is an open neighborhood of $m$. 
So, if  $(z_{\alpha})$ is a net in $S$ converging to $m$ then there is $\alpha (k)$ 
such that the tail $(z_{\alpha})_{ \alpha \geq \alpha (k) }$ is completely contained in 
$\{ z\in \disc : |b_{k}(z)| < a/2 \} \cap S = S_{k}$.  
Therefore (\ref{ely}) implies that   
$\lim_{\alpha} g (L_{z_{\alpha}} (z) )  = f(z)$.

By (3), $\sum_{k,n \geq 1} (1-|z_{k,n}|)  <  \sum_{k\geq 1} 1/2^{k} = 1$, 
and  consequently the Blaschke product 
$b= \prod_{k\geq 1}  b_{k}$ converges. 
Furthermore, if we write $\beta = \prod_{k\geq 1} \sigmak$, condition (5) tells us that 
\[ |b(z)| \geq \beta \ \ \mbox{ when }\ \ 
z\not \in \bigcup_{k\geq 1}  \bigcup_{n\geq 1}  \hypdi ( z_{k,n}  , s_{k})  = \Om  . 
\] 
That is, $V = \{ z\in \disc : |b(z)| < \beta \}  \subset \Om$. In addition, since each 
$b_{k}$ vanishes on $m$ (because $m\in \ov{S}_{k}$ for every $k\geq 1$) then
$b$ vanishes on $m$ with infinite multiplicity. So, $b\equiv 0$ on $P(m)$. 

Thus, $g\in \papa (\Om ) \subset \papa (V)$. 
By Lemma \refth{useful} then there are
$0< \gamma = \gamma (\beta ) < \beta$, a constant $C>0$ and $F\in \papa$ such that 
\begin{equation}    \laba{prelem} 
|F(z) -g(z)|  \leq   C |b(z)|   \ \ \mbox{ when }\ \  |b(z)| < \gamma  . 
\end{equation} 
%
Let $(z_{\alpha})$ be any net in $S$ that tends to $m$ and let $z\in \disc$. 
Then $L_{z_{\alpha}} (z) \rr L_{m}(z) \in P(m)$, and since $b\equiv 0$ on $P(m)$ then 
$b\circ L_{z_{\alpha}} (z)  \rr b\circ L_{m}(z) =0$. 
So, there is $\alpha_{0}$ (depending on $z$) such that 
$L_{z_{\alpha}}(z) \in \{ |b| < \gamma \}$   for every $\alpha \geq \alpha_{0}$. 
Thus, by (\ref{prelem}) 
\[   |F\circ L_{z_{\alpha}} (z)  -  g\circ L_{z_{\alpha}} (z)|  \leq |b\circ L_{z_{\alpha}} (z)| 
\ \ \mbox{ for }\ \ \alpha \geq \alpha_{0} , 
\]
where the last term tends to zero when we take limit in $\alpha$. 
Henceforth  
\[
F\circ L_{m} (z)  =  \lim_{\alpha} F\circ L_{z_{\alpha}}(z)  
			    =  \lim_{\alpha} g\circ L_{z_{\alpha}}(z)  = f(z)  , 
\] 
and we are done.    \diaend  \\ 

\noindent 
The proof above shows that if $y\in \bigcap_{k\geq 1} \ov{S_{k}}$ 
is any point and $F\in \papa$ is the function constructed in the last step, then 
$F\circ L_{y} (z) = f(z)$. 
This does not mean that ${\cal L}_{y} (0) = {\cal L}_{m} (0)$, 
because the chain of interpolating sequences constructed depends 
on the function $f$.  

\section{Examples} 
\equnew 

A point $m\in \chu \setminus \disc$ is called {\em oricycular} if it is in the closure 
of a region limited by two circles in $\disc$ that are tangent to $\partial \disc$ at the 
same point. Every oricycular point is in $\calg$ and it is in the closure 
of some tangent circle to $\partial \disc$ (see \cite[pp.$\,$107-108]{hof} ). 
We are going to search for the possible fibers that meet ${\cal L}_{m}(0)$ when $m$ is 
a nontangential or an oricycular point, and we shall determine $\lambda (m,x)$ 
for all possible $x\in {\cal L}_{m}(0)$. 

Every point in $\chu \setminus \disc$ has the form 
$\gamma m$, where $\pi (m) =1$ and 
$\gamma \in {\Bbb C}$ has modulus 1.  
Since $\lambda (m,  x) = \lambda (\gamma m, \gamma x )$ 
for every $x\in \chu$, and 
by Lemma \refth{mul}  
${\cal L}_{\gamma  m} (0)  = \gamma {\cal L}_{m}(0)$, 
then  
there is no loss of generality by considering $\pi (m) =1$. 

Let $b$ be an interpolating Blaschke product with zero sequence $\{ z_{k} \}$ 
such that $b(m) =0$. 
If  the point  $\om \in \disc$ is a zero of  $b\circ L_{m}$ then  
there is a subsequence $\{ z_{k_{j}} \}$ of  $\{ z_{k} \}$ such that 
$b\circ L_{ z_{k_{j}} } (\om ) \rr 0$. 
By Lemma \refth{dist} then   
\[
  \lim_{j}  \rho (   \om  ,   \{   L_{-z_{k_{j}}} (  z_{k} )   \}_{k\geq 1}  ) =  
\lim_{j}   \rho (  L_{z_{k_{j}}}(\om )  ,   \{ z_{k} \}  ) 
\leq 
c(\delta(b))^{-1}    \lim_j  |b(   L_{z_{k_j}} (\om)  )| = 0 .   
\]   
Consequently, $\om$ is an accumulation point of 
\[  A_{N} =   \{  L_{ -z_{n} } (z_{k}) :   k \geq 1 ,  \      n \geq N  \} 
\]
for every positive integer $N$. 

If $m$ is a nontangential point we can assume that 
there is some fixed $-\pi /2 < \theta < \pi /2$, such that $\{ z_{n} \}$ lies in the 
straight segment 
\[    S = \{ z\in \disc : 1-z = re^{i\theta} ,  \  r > 0   \}   . 
\] 
A straightforward calculation shows that the closure of $S$ in ${\Bbb C}$ 
meets $\partial \disc$ when $r=0$ and $r=2\cos \theta$. 
Therefore,  
$S = ( 1-2\cos \theta \, e^{i\theta} , 1 )$ and 
$z_{n} = 1-r_{n}\eiti$, with $0< r_{n} < 2\cos \theta$. 
We can also assume that $z_n \rr 1$. 
The conformal map $L_{-z_{n}}$ sends $S$ into a circular segment 
$C_{n} \cap \disc$, where $C_{n}$ is the circle that pass through the points 
$L_{ -z_{n} } ( z_{n})  = 0$,             
$L_{ -z_{n} } (1)          =  e^{i2\theta}$  and  
\begin{equation}      \laba{veryla}
L_{ -z_{n} } (1-2\cos \theta \, \eiti )    =   
\frac{ (r_{n} -2\cos \theta ) \eiti }{ 
2\cos \theta \, \eiti + r_{n} (e^{-i\theta} - 2\cos \theta )} .     
\end{equation}  
We are including here the extreme case when   
$C_{n}$ is a straight line (i.e., $\theta = 0$). 
Since $r_{n} \rr 0$, taking  limits  in 
(\ref{veryla}) we see that the limit curve of  $C_{n}$ 
is the circular segment $C\cap \disc$, where $C$ is the circle that pass through 
$0$,  $e^{i2\theta}$ and $-1$. 
Therefore the zero sequence of $b\circ L_{m}$ lies in $C\cap \disc$, and 
since $b\circ L_{m}$  vanishes on ${\cal L}_{m}(0)$, only the fibers of $-1$ and 
$e^{i2\theta}$ can have points of  ${\cal L}_{m}(0)$. 
 
Clearly, if $x\in {\cal L}_{m}(0)$ is in the fiber of $e^{i2\theta}$ then 
$\lambda (m,x ) = e^{i2\theta}$. 
Suppose that  $x\in {\cal L}_{m}(0)$ is in the fiber of $-1$. 
Since the straight segment 
$-1 + R e^{-i\theta}$,  with  $0< R< 2\cos \theta$, is tangent to $C$ at the point 
$-1$ then $x$ is a nontangential point lying in the closure of this segment. 
So, $\lambda (m,x) = e^{-i2\theta}$ by Section 1. 

If $m$ is a oricycular point (with $\pi (m) =1$) a similar but easier analysis shows that 
${\cal L}_{m}(0)$ lies in the closure of  $C\cap \disc$, where $C$ is the tangent circle to 
$\partial \disc$ that pass through $0$ and $-1$.  
Therefore ${\cal L}_{m}(0)$ only can meet the fiber of $-1$, and indeed it does unless 
$P(m)$ is a homeomorphic disk. 
Since $C$ is tangent to $\partial \disc$ then every $x\in {\cal L}_{m}(0)$ with 
$\pi (x) = -1$ is a tangential point, which by Section 1 yields 
$\lambda (m ,x ) = -1$. 
%
%
%
%

\newcommand{\foo}{\footnotesize}

 \noindent Daniel Su\'{a}rez\\
 Departamento de Matem\'{a}tica \\
 Facultad de Cs.\@ Exactas y Naturales \\
 UBA, Pab. I, Ciudad Universitaria \\
 (1428) N\'{u}\~{n}ez, Capital Federal \\
 Argentina\\
\vspace{0.5mm}
\noindent
$\! \!${\foo dsuarez@dm.uba.ar}

\end{document}